\documentclass[aos,MSNbibl,nameyear,seceqn,dvips]{arximspdf}

\doi{10.1214/11-AOS900}
\volume{39}
\issue{4}
\pubyear{2011}
\firstpage{2103}
\lastpage{2130}

\makeatletter

\newcommand{\iint}{\int\!\!\int}

\newcommand{\xrightarrow}[1]{\stackrel{#1}{\longrightarrow}}

\newcommand{\C}{{\mathbb C}}
\newcommand{\N}{{\mathbb N}}
\newcommand{\Z}{{\mathbb Z}}
\newcommand{\imag}{\mathrm i}
\newcommand{\ve}{\varepsilon}

\newtheorem{lemma}{Lemma}[section]
\newtheorem{theorem}{Theorem}[section]
\newtheorem{corollary}{Corollary}[section]

\newproclaim{remark}{Remark}[section]
\newproclaim{example}{Example}[section]

\makeatother

\begin{document}
\begin{frontmatter}

\title{On the range of validity of the autoregressive sieve bootstrap}
\runtitle{On the validity of the autoregressive sieve bootstrap}

\begin{aug}
\author[A]{\fnms{Jens-Peter} \snm{Kreiss}},
\author[B]{\fnms{Efstathios} \snm{Paparoditis}}
and
\author[C]{\fnms{Dimitris~N.}~\snm{Politis}\corref{}\ead[label=e1]{dpolitis@ucsd.edu}}
\runauthor{J.-P. Kreiss, E. Paparoditis and D. N. Politis}
\affiliation{Technische Universit\"at Braunschweig, University of
Cyprus and University of California, San Diego}
\address[A]{J.-P. Kreiss\\
Institut f\"ur Mathematische Stochastik\\
Technische Universit\"at Braunschweig\\
Pockelsstrasse 14\\
D-38106 Braunschweig\\
Germany}
\address[B]{E. Paparoditis\\
Department of Mathematics\\
\quad and Statistics\\
University of Cyprus\\
1678 Nicosia\\
Cyprus}
\address[C]{D. N. Politis\\
Department of Mathematics\\
University of California, San Diego\\
La Jolla, California 92093--0112\\
USA\\
\printead{e1}}
\end{aug}

\received{\smonth{12} \syear{2010}}
\revised{\smonth{4} \syear{2011}}

%
\begin{abstract}
We explore the limits of the autoregressive (AR) sieve bootstrap, and
show that its applicability extends well beyond the realm of linear
time series as has been previously thought. In particular, for
appropriate statistics, the AR-sieve bootstrap is valid for
stationary processes possessing a general Wold-type autoregressive
representation with respect to a white noise; in essence, this includes
all stationary, purely nondeterministic processes, whose spectral
density is everywhere positive. Our main theorem provides a simple and
effective tool in assessing whether the AR-sieve bootstrap is asymptotically
valid in any given situation. In effect, the large-sample distribution of
the statistic in question must only depend on the first and second order moments
of the process; prominent examples include the sample mean and the spectral density.
As a counterexample, we show how the AR-sieve bootstrap is not always valid for
the sample autocovariance even when the underlying process is linear.
\end{abstract}

%

\vspace*{12pt}
\end{frontmatter}

\section{Introduction}

Due to the different possible dependence structures that may occur in
time series analysis, several bootstrap procedures have been proposed
to infer properties of a statistic of interest. Validity of the
different bootstrap procedures depends
on the probabilistic structure of the underlying stochastic process
${\mathbf X}= (X_t \dvtx t\in\mathbb{Z} ) $
and/or on the particular statistic considered. Bootstrap schemes for
time series rank from those imposing more parametric type assumptions on
the underlying stochastic process class to those accounting only for
some kind of mixing or weak dependence assumptions. For an overview see
\citet{Buhlmann02}, \citet{Lahiri2003}, \citet
{Politis03} and \citet{PapPol09}.

A common assumption is that $ {\mathbf X}$ is a \textit{linear time series},
that is, that
%
%
\begin{equation}
\label{MAinfinity}
X_t= \sum_{j=-\infty}^{\infty} b_j e_{t-j} ,\qquad t\in\mathbb{Z} ,
\end{equation}
with respect to independent, identically distributed (i.i.d.) random
variables $(e_t)$---often assumed to have mean zero and finite fourth
order moments---and for absolutely summable coefficients $(b_j)$;
this is not to be confused with\vadjust{\goodbreak} the Wold representation
with respect to white noise, that is, uncorrelated, errors that all
stationary, purely nondeterministic processes possess. If $b_j=0$ for
all $j<0$, then the linear process is called \textit{causal}.

Stationary autoregressive (AR) processes of order $p$
are members of the linear class (\ref{MAinfinity}) provided the
autoregression is defined
on the basis of i.i.d. errors. Model-based bootstrapping in the AR($p$)
case was among the
first bootstrap proposals for time series; see, for example, \citet
{Freedman1984}.
The extension to the AR($\infty$) case was inevitable;
this refers to the situation where the strictly stationary process
$ X_t$ has the following
linear infinite order autoregressive representation
%
%
\begin{equation} \label{eqarinfty}
X_t = \sum_{j=1}^{\infty} \pi_j X_{t-j} + e_t ,\qquad t\in\mathbb{Z} ,
\end{equation}
with respect to i.i.d. errors $e_t$ having mean zero,
variance $ 0< E( e_t^2) = \sigma^2_e $ and $ E(e_t^4)< \infty$;
here the coefficients $ \pi_j$ are assumed absolutely summable and $ \pi
(z)=1 - \sum_{j=1}^{\infty}\pi_j z^j \neq0$ for $ |z| = 1$.
The two representations, (\ref{MAinfinity}) and~(\ref{eqarinfty}) are
related; in fact,
the class~(\ref{eqarinfty}) is a subset of the linear class~(\ref{MAinfinity}).
Furthermore, it can be shown that the linear AR($\infty$) process (\ref
{eqarinfty}) is causal
if and only if $ \pi(z)=1 - \sum_{j=1}^{\infty}\pi_j z^j \neq0$ for $
|z| \leq1$.

There is already a large body of literature
dealing with applications and properties of the AR-sieve bootstrap.
Kreiss (\citeyear{Kreiss88}, \citeyear{Kreiss92}) established validity
of this bootstrap scheme for
different statistics including
autocovariances and autocorrelations. \citet{PapStreit91}
established asymptotic validity of the AR-sieve
bootstrap to infer properties of high order autocorrelations, and
\citet{Pap96} established its
validity in a multivariate time series context. The aforementioned results
required an exponential decay of the AR coefficients $ \pi_j$ as $ j
\rightarrow\infty$; \citet{Buhlmann97} extended the class of
AR($\infty$) processes for which the AR-sieve bootstrap works by allowing
a polynomially decay of the $\pi_j$
coefficients. Furthermore, \citet{BickelBuhlmann99} introduced a
mixing concept appropriate for
investigating properties of the AR-sieve bootstrap which is related to
the weak dependence concept of \citet{DoukhanLouhichi99},
while \citet{ChoiHall00} focused on properties of the AR-sieve
bootstrap-based confidence intervals.

A basic assumption in the current literature of the AR-sieve bootstrap
is that ${\mathbf X}$ is a linear
AR($\infty$) process, that is, $ X_t$ is generated by
(\ref{eqarinfty}) with $(e_t)$ being an i.i.d. process.
One exception is the case of the sample mean $ \overline
{X}_n=n^{-1}\sum_{t=1}^{n}X_t$, where
\citet{Buhlmann97} proved
validity of the AR-sieve bootstrap also for the case where the
assumption of i.i.d. errors in (\ref{eqarinfty})
can be relaxed to that of martingale differences, that is,
$E(e_t|{\mathcal E}_{t-1}) = 0$ and $E(e_t^2|{\mathcal E}_{t-1}) =\sigma^2_e$
with ${\mathcal E}_{t-1}=\sigma(\{ e_s \dvtx s \leq t-1\})$ the $\sigma$-algebra generated by the random
variables $\{ e_{t-1}, e_{t-2},\ldots\}$. Notice that the process (\ref
{eqarinfty}) with
innovations forming a martingale difference sequence is in some sense
not ``very far'' from the linear process
(\ref{eqarinfty}) with i.i.d. errors.
In fact, some authors call the set-up of model (\ref{MAinfinity})
with martingale difference errors
``weak linearity,'' and the same would hold regarding (\ref
{eqarinfty}); see, for example, \citet{KokoPoli2011}.

To elaborate, for a causal linear process
the general $L_2$-optimal
predictor of $X_{t+k}$ based on its past $ X_t, X_{t-1}, \ldots,$
namely the conditional expectation
$ E(X_{t+k}|X_s, s \leq t)$ of $ X_{t+k}$, is
identical to the best \textit{linear} predictor~${\mathcal P}_{{\mathcal M}_{t}}(X_{t+k})$;
here $k$ is assumed positive, and $ {\mathcal P}_{C}$ denotes
orthogonal projection onto the set $C$ and
${\mathcal M}_{s}=\overline{\mathrm{span}}\{X_j \dvtx j \leq s\}$, that is,
the closed linear span generated by the random variables
$\{X_j \dvtx j \leq s\}$. The property of linearity of the optimal
predictor is shared by causal processes
that are only weakly linear.
Recently, under the assumption of weak linearity with (\ref
{eqarinfty}), \citet{Poskitt07} claimed
validity of the AR-sieve bootstrap
for a much wider class of statistics that are defined as smooth
functions of means.
However, this claim does not seem to be correct in general. In
particular, our Example~\ref{bootautocov}
of Section \ref{bootfctlinstat} contradicts Theorem 2 of \citet
{Poskitt07}; see Remark~\ref{rePoskitt} in what follows.

The aim of the present paper is to explore the limits of the AR-sieve
bootstrap, and to give a definitive answer to the question concerning
for which classes of statistics, and for which dependence structures,
is the AR-sieve bootstrap asymptotically valid. Moreover, we also
address the question what the AR-sieve bootstrap really does when it is
applied to data stemming from a stationary process not fulfilling
strict regularity assumptions such as linearity or weak linearity. In
order to do this, we examine in detail in Section~\ref{generalAR}
processes possessing a so-called general autoregressive representation
with respect to white noise errors; these form a much wider class of
processes than the linear AR($\infty$) class described by
(\ref{eqarinfty}).

Our theoretical results in Section \ref{bootfctlinstat} provide an
effective and simple tool for gauging consistency of the AR-sieve
bootstrap. They imply that for certain classes of statistics the range
of the validity of the AR-sieve bootstrap goes far beyond that of the
linear class (\ref {MAinfinity}). On the other hand, for other classes
of statistics, like for instance autocorrelations, validity of the
AR-sieve bootstrap is restricted to the linear process class (\ref
{MAinfinity}), while for statistics like autocovariances, the
\mbox{AR-sieve} bootstrap is only valid for the linear AR($\infty$)
class (\ref{eqarinfty}). But even in the case of the linear
autoregression (\ref{eqarinfty}) with infinite order, the theory
developed in this paper provides a further generalization of existing
results since it establishes validity of this bootstrap procedure under
weaker assumptions on the summability of the coefficients $\pi_j$, thus
relaxing previous assumptions referring to exponential or polynomial
decay of these coefficients.

The remaining of the paper is organized as follows. Section \ref
{generalAR} develops the background concerning the Wold-type infinite
order AR representation that is required to study the AR-sieve
bootstrap, and states the necessary assumptions to be imposed on the
underlying process class and on the parameters of the bootstrap
procedure. Section \ref{bootfctlinstat} presents our main result and
discusses its implications by means of several examples. Proofs and
technical details are deferred to the \hyperref[proofs]{Appendix}.

\section{The AR-sieve bootstrap and general autoregressive representation}
\label{generalAR}

Here, and throughout the paper, we
assume that we have observations~$X_1,\ldots,\allowbreak X_n$ stemming from a
strictly stationary process $ {\mathbf X}$.
Let $T_n=T_n(X_1,\ldots,X_n)$ be an estimator of some unknown
parameter $\theta$ of the underlying stochastic
process $ {\mathbf X}$. Suppose that for some appropriately increasing
sequence of real numbers
$ \{c_n \dvtx n \in\N\}$ the distribution $ \mathcal{L}_n=\mathcal
{L}(c_n(T_n - \theta))$
has a nondegenerated limit.
The AR-sieve bootstrap proposal to estimate the distribution~$\mathcal{L}_n$ goes as follows:
\begin{longlist}[Step 3:]
\item[Step 1:] Select an order $p=p(n) \in\mathbb{N}$, $ p \ll n$,
and fit a $p$th order autoregres\-sive model to $ X_1, X_2, \ldots,
X_n$. Denote by
$ \widehat{a}(p)=(\widehat{a}_j(p), j=1,2, \ldots, p)$,\vspace*{1pt} the Yule--Walker
autoregressive parameter estimators, that is,
$\widehat{a}(p)=\widehat{\Gamma}(p)^{-1}\widehat{\gamma}_p$ where for $
0 \leq h \leq p$,
\[
\widehat{\gamma}_X(h)=n^{-1}\sum_{t=1}^{n-|h|}(X_t-\overline{X}_n)
\bigl(X_{t+|h|}-\overline{X}_n\bigr),
\]
$\overline{X}_n=n^{-1}\sum_{t=1}^{n}X_t$, $ \widehat{\Gamma
}(p)=(\widehat{\gamma} _X(r-s))_{r,s=1,2, \ldots, p}$ and
$ \widehat{\gamma}_p=( \widehat{\gamma}_X(1), \widehat{\gamma}_X(2),
\ldots,\break \widehat{\gamma}_X(p))^{\prime}$.

\item[Step 2:] Let
$ \widetilde{\varepsilon}_t(p) =X_t
-\sum_{j=1}^{p}\widehat{a}_j(p)X_{t-j}$,
$ t=p+1, p+2, \ldots, n$, be the residuals of the autoregressive fit
and denote by $ \widehat{F}_n$ the empirical distribution function of
the centered residuals
$ \widehat{\varepsilon}_t(p)=\widetilde{\varepsilon}_t(p)-\overline
{\varepsilon}$, where $ \overline{\varepsilon}=(n-p)^{-1}\sum
_{t=p+1}^{n}\widetilde{\varepsilon}_t(p)$.
Let $( X_1^\ast, X_2^\ast, \ldots, X_n^\ast)$ be a set of observations
from the time series $ {\mathbf X}^{\ast}=\{X_t^\ast\dvtx t \in\mathbb
{Z}\}$ where
$ X^{\ast}_t = \sum_{t=1}^{p}\widehat{a}_j(p)X^{\ast}_{t-j} + e^\ast_t$
and the $ e^\ast_t$'s are independent random variables
having identical distribution $ \widehat{F}_n$.

\item[Step 3:] Let $ T_n^\ast=T_n(X_1^\ast,X_2^\ast, \ldots,
X_n^\ast)$ be the same
estimator as $ T_n$ based on the pseudo-time series $ X_1^\ast, X_2^\ast
, \ldots, X_n^\ast$
and $ \theta^\ast$ the analogue of $ \theta$ associated with the
bootstrap process
${\mathbf X}^\ast$. The AR-sieve bootstrap approximation of~$\mathcal
{L}_n$ is then given by
$ \mathcal{L}^\ast_n= \mathcal{L}^\ast(c_n (T_n^\ast-\theta^\ast))$.
\end{longlist}

In the above (and in what follows), $\mathcal{L}^\ast, E^\ast, \ldots$
will denote
probability law, expectation, etc. in the bootstrap world (conditional
on the data $X_1,\ldots,X_n$).\vadjust{\goodbreak}

Note that the use of Yule--Walker estimators in Step 1 is essential and
guarantees---among other things---that the complex polynomial
$ \widehat{A}_p(z) = 1 - \sum_{j=1}^p\widehat a_j(p)z^j$ has no roots
on or within the unit disc
$ \{z\in\C\dvtx|z| \leq1\}$, see the discussion before (\ref
{zeroesAhat}), that is, the bootstrap process
$ {\mathbf X}^\ast$ always is a~sta\-tionary and causal autoregressive process.

The question considered in this paper is when can the bootstrap
distribution $ \mathcal{L}_n^{\ast}$ correctly
approximate the distribution $ \mathcal{L}_n$ of interest, and moreover what
the AR-sieve bootstrap does if the latter is not the case. To this
end, let us first discuss a general autoregressive representation of
stationary processes.

Recall that by the well-known Wold representation, every purely
nondeterministic, stationary and zero-mean stochastic process
$ {\mathbf X}=\{X_t \dvtx t \in\mathbb{Z}\}$ can be expressed as
%
%
\begin{equation} \label{eqwold}
X_t = \sum_{j=1}^{\infty} b_j u_{t-j} + u_t,
\end{equation}
where $\sum_{j=1}^{\infty}b_j^2 < \infty$ and $ u_t =X_t-{\mathcal
P}_{{\mathcal M}_{t-1}}(X_t)$ is a zero mean,
white noise ``innovation'' process with finite variance $0< \sigma
^2_u=E(u_t^2) < \infty$;
recall that ${\mathcal M}_{s}=\overline{\mathrm{span}}\{X_j \dvtx j
\leq s\}$.

Less known is that for all purely nondeterministic, stationary and
zero-mean time series
unique autoregressive coefficients
$(a_k\dvtx k\in\N)$ exist that only depend on the autocovariance function
of the time series $(X_t)$,
such that for any $ n \in\N$,
%
%
\begin{equation}
\label{arcoeff}
\mathcal{P}_{\mathcal{M}_{t-1}}(X_t) = \sum_{k=1}^na_k
X_{t-k}+e_{t,n} ,\qquad t\in\mathbb{Z} ,
\end{equation}
where $(e_{t,n} \dvtx t\in\mathbb{Z})$ is stationary and $e_{t,n}\in
\overline{\mathrm{span}}\{ X_s \dvtx s\le t-n-1\}$.

Under the additional assumption that the coefficients $ (a_k, k \in\N
)$ are absolute summable, that is,
$\sum_{k=1}^{\infty} |a_k| <\infty$, one then obtains an
autoregressive, Wold-type
representation of the underlying process given by
%
%
\begin{equation}
\label{woldar}
X_t = \sum_{k=1}^{\infty} a_k X_{t-k} + \ve_t ,\qquad t\in\mathbb{Z} .
\end{equation}
Here again $(\ve_t\dvtx t\in\mathbb{Z})$ denotes a white noise, that is,
uncorrelated,
process with finite variance $\sigma_{\ve}^2=E \ve_t^2$ which fulfills
%
%
\begin{equation}
\label{sigmasquare}
\sigma_{\ve}^2 = \gamma_X(0) - \sum_{k=1}^{\infty} a_k \gamma_X(k),
\end{equation}
where $\gamma_X(\cdot)$ denotes the autocovariance function of $
{\mathbf X}$.

Under the absolute summability assumption on the autoregressive
coefficients $(a_k)$---conditions for which will be
given in Lemma \ref{lemma1} in the sequel---we have that\vadjust{\goodbreak}
$ X_t-{\mathcal P}_{{\mathcal M}_{t-1}}(X_t)
=X_t-\sum_{k=1}^{\infty} a_k X_{t-k} $; this implies that the
white noise process $(u_t)$ appearing in
(\ref{eqwold})
coincides with the white noise process $(\ve_t)$ in (\ref{woldar}).
Notice that this does not mean that if we have an arbitrary one sided
moving average representation of a time series $(X_t)$, even with
summable coefficients, that this
moving average representation is the Wold representation of the
process; see Remark \ref{remark1} for an example.
Furthermore, let $ f_X$ be the spectral density of $ {\mathbf X}$, that is,
\[
f_X(\lambda)=(2\pi)^{-1}\sum_{h\in\Z}\gamma(h)\exp\{-\imag\lambda
h\},\qquad
\lambda\in[-\pi,\pi].
\]
Then, from (\ref{woldar}) one immediately obtains that
%
%
\begin{equation}
\label{Azeroesfirst}
\Biggl\vert1-\sum_{k=1}^{\infty}a_k e^{-\imag k\lambda} \Biggr\vert^2
\cdot f_X(\lambda) =
\frac{\sigma_{\ve}^2}{2\pi} ,\qquad \lambda\in[-\pi, \pi],
\end{equation}
which implies that for strictly positive spectral densities $f_X$
the power series $A(z):=1-\sum_{k=1}^{\infty}a_k z^k$ has no zeroes
with $|z|=1$.
For more details of the autoregressive Wold representation (\ref{woldar})
see \citet{Pourahmadi01}, Lemma~6.4(b), (6.10) and (6.12).
It is worth mentioning that in the historical evolution of Wold
decompositions the autoregressive variant preceded the moving average one.

\begin{remark}
\label{remark1}
If we consider a purely nondeterministic and stationary time series
possessing a standard one-sided
moving average representation, and if we additionally
assume that the spectral density is bounded away from zero and that the
moving average coefficients $b_j$
are absolutely summable, then this would imply that the polynomial
$B(z)=1+\sum_{j=1}^{\infty}b_jz^j$ has no zeroes with magnitude equal
to one.
There may of course exist zeroes within the unit disk. But since the
closed unit disk is compact and $B(z)$ represents a holomorphic
function there could exist only finitely many zeroes with magnitude
less than one. Following the technique described in Kreiss and Neuhaus
[(\citeyear{KreissNeuhaus06}), Section 7.13] one may switch to another moving
average model for which the polynomial has no zeroes within the unit
disk. This procedure definitely changes the white noise process; for
example, if the white noise process in the
assumed moving average representation consists of independent random
variables, this desirable feature typically is lost when switching to
the moving average model with all zeroes within the unit disk removed.
In fact, only the property of uncorrelatedness is
preserved.
The modified moving average process allows then for an autoregressive
representation of infinite order and this process, because of the
uniqueness of the autoregressive representation, coincides with the one
in (\ref{woldar}).

The following simple example, taken from Brockwell and Davis
[(\citeyear{BrockwellDavis91}),
Example 3.5.2] illustrates these points.
Based on i.i.d. random variables $(e_t)$ with mean zero and finite and
nonvanishing variance $\sigma^2_e$, construct\vadjust{\goodbreak}
the simple MA(1)-process
%
%
\begin{equation}
\label{MA1}
X_t= e_t-2e_{t-1},\qquad t\in\mathbb{Z} .
\end{equation}
This MA(1)-model is not invertible to an autoregressive process.
However, a~general autoregressive representation as described above exists.
In order to obtain this representation
denote by $L$ the usual lag-operator and consider
$B(L):= 1-2L$ as well as $\widetilde B(L):=1-0.5 L$. Of course
%
%
\begin{equation}
X_t= \widetilde B(L) \frac{B(L)}{\widetilde B(L)} e_t .
\end{equation}
Since $|B(e^{-\imag\lambda})|^2 / |\widetilde B(e^{-\imag\lambda
})|^2 = 4$, we obtain that
\[
\ve_t := \frac{B(L)}{\widetilde B(L)} e_t = e_t -\frac{3}{2}\sum
_{j=1}^{\infty}\biggl(\frac{1}{2}\biggr)^{j-1} e_{t-j}.
\]
Again $(\varepsilon_{t}) $ is a (uncorrelated) white noise process with variance
$\sigma_{\ve}^2=4 \sigma_e^2$. Moreover, we have
%
%
\begin{equation}
\label{alternativeMA1}
X_t=\ve_t -0.5 \ve_{t-1} = -\sum_{j=1}^{\infty} 0.5^j X_{t-j} + \ve_t.
\end{equation}
Obviously $\ve_t= X_t-\mathcal{P}_{\mathcal{M}_{t-1}}(X_t)$ which
means that (\ref{alternativeMA1}) and not (\ref{MA1}) is the Wold
representation of the time series $(X_t)$. This also means that the
modified moving average (or Wold) representation of the
process $(X_t)$ possesses only uncorrelated innovations $(\ve_t)$
instead of independent innovations~$(e_t)$. But the representation (\ref
{alternativeMA1}) with uncorrelated innovations has the
advantage that it
indeed possesses an autoregressive representation of infinite order. Of course,
via the described modification, we do not change any property of the
process $(X_t)$.
But, and this is essential, the modification leading to the general
AR($\infty$)-representation typically destroys a existing independence
property of the white noise in a former moving average representation.

To elaborate, the problem of understanding the stochastic properties
of the innovation process in linear time series has been thoroughly
investigated in the literature.
Breidt and Davis (\citeyear{BreidtDavis91}) showed that time reversibility of a~linear
process is equivalent to the fact that the i.i.d. innovations $ e_t$
are Gaussian
and used this result to derive for a class of linear processes
uniqueness of moving average representations with i.i.d. non-Gaussian
innovations and to discuss the stochastic properties of the innovation
process appearing in alternative moving average representations for the
same process class.
\citet{Breidtetal95} used such results to initialize autoregressive
processes in Monte Carlo generation of conditional sample paths running
autoregressive processes backward in time and \citet
{Andrewsetal07} for
estimation problems for all-pass time series models. Properties of the
innovation process
in non-Gaussian, noninvertible time series have been also discussed in
Lii and Rosenblatt (\citeyear{LiiRosen82,LiiRosen96}).

As we have seen the variances of $e_t$ and $\ve_t$ do not coincide and
the same is true for the fourth order cumulant
$E (e_t^4)/\sigma_e^4 -3$ which will be of some importance later.
Using the fact that $\ve_t$ is defined via a linear transformation on
the i.i.d. sequence $(e_t)$
we obtain by straightforward computation
%
%
\begin{equation}
\label{cumulanteps}
\frac{E (\ve_t^4)}{(E (\ve_t^2))^2}-3 = \frac{2}{5} \frac{E
(e_1^4)}{\sigma_e^4}-\frac{6}{5} ,
\end{equation}
which only equals $E (e_1^4)/\sigma_e^4 -3$ in case the latter
quantity is equal to $0$, for example, when the $e_t$ are normally distributed.
The normally distributed case always leads to the fact that
uncorrelatedness and independence are equivalent, thus implying that the
white noise process in the general
autoregressive representation always consists of independent and
normally distributed random variables which leads for the autoregressive
sieve bootstrap in some cases to a considerable simplification as we
will see later.
\end{remark}

In order to get conditions which ensure the absolute summability of the
autoregressive coefficients $(a_k, k \in\N)$, one can go back to
an important paper
by \citet{Baxter62}. Informally speaking it is the smoothness of
the spectral density $f_X$ which ensures summability of these
coefficients. To be more precise, we have the following result.
\begin{lemma}
\label{lemma1}
\textup{(i)} If $f_X$ is strictly positive and continuous and if
\[
\sum
_{h=0}^{\infty} h^r |\gamma_X(h)| < \infty
\]
for some $r \ge0$, then
%
%
\begin{equation}
\sum_{h=0}^{\infty} h^r |a_h| < \infty.
\end{equation}

\textup{(ii)} If $f_X$ is strictly positive and possesses $k\ge2$
derivatives, then
%
%
\begin{equation}
\sum_{h=0}^{\infty} h^r |\gamma_X(h)| < \infty\qquad \forall r
< k-1.
\end{equation}
\end{lemma}
\begin{pf}
Cf. \citet{Baxter62}, pages 140 and 142.
\end{pf}

The uniquely determined autoregressive coefficients $(a_k)$ are closely
related to
the coefficients of an optimal (in the mean square sense)
autoregressive fit of order~$p$, or equivalently, to
prediction coefficients based on the finite past. To be precise, denote
the minimizers of
%
%
\begin{equation}
E \Biggl( X_t - \sum_{r=1}^p c_r X_{t-r}\Biggr) ^2
\end{equation}
by $a_1(p),\ldots, a_p(p)$, which of course are solutions of the
following Yule--Walker linear equations:
%
%
\begin{eqnarray}
\label{ARplineareq}
\pmatrix{
\gamma_X(0) & \cdots& \gamma_X(p-1) \cr
\vdots& \ddots& \vdots\cr
\gamma_X(p-1) & \cdots& \gamma_X(0)}
\pmatrix{
c_1 \cr
\vdots\cr
c_p}
&=&
\pmatrix{
\gamma_X(1) \cr
\vdots\cr
\gamma_X(p)}
.
\end{eqnarray}
Recall from Brockwell and Davis [(\citeyear{BrockwellDavis91}),
Proposition 5.1.1] that the
covariance matrix $\Gamma(p)$ on the left-hand side is for all
$p$ invertible provided $\gamma_X(0)>0$ and $\gamma_X(h) \to0$ as
$h\to\infty$.

Now by slight modifications of \citet{Baxter62}, Theorem 2.2 [cf.
also Pourahmadi (\citeyear{Pourahmadi01}),
Theorem 7.22], we obtain the following helpful result relating the
coefficients $ a_k(p)$ of the $p$th order autoregressive
fit to the $ (a_k)$ of the general autoregressive representation.
\begin{lemma}
\label{lemma2}
Assume that $f_X$ is strictly positive and continuous and that $\sum
_{h=0}^{\infty} (1+h)^r |\gamma_X(h)| < \infty$ for some $r\ge0$.
Then there exists $p_o\in\mathbb{N}$ and $C>0$ (both depending on
$f_X$ only) such that for all $p\ge p_o$,
%
%
\begin{equation}
\sum_{k=0}^{p} (1+k)^r |a_k(p)-a_k| \le C\cdot\sum_{k=p+1}^{\infty}
(1+k)^r |a_k|
\end{equation}
as well as
%
%
\begin{equation}
\sum_{k=1}^{\infty} (1+k)^r |a_k| < \infty.
\end{equation}
This means that we typically can achieve a polynomial rate of
convergence of $a_k(p)$ toward $a_k$.
\end{lemma}


As already mentioned, 
$\gamma_X(0)>0$ and
$\gamma_X(h)\to0$ as $h\to\infty$ ensure nonsingularity of all
autocovariance matrices appearing in the left-hand side
of~(\ref{ARplineareq}). Since these matrices
are positive semidefinite this means that under these conditions
$\Gamma(p)$ actually is positive definite.
This in turn with Kreiss and Neuhaus [(\citeyear{KreissNeuhaus06}), Section 8.7]
implies
that the polynomial
$A_p(z)=1-\sum_{k=1}^p a_k(p)z^k$ has no zeroes in the closed unit
disk. We can even prove a slightly stronger result.
\begin{lemma}
\label{lemma3}
Assume that $f_X$ is strictly positive and continuous, that $\sum
_{h=0}^{\infty} |\gamma_X(h)| < \infty$
and $\gamma_X(0)>0$.\vadjust{\goodbreak}
Then there exists $\delta>0$ and $p_o\in\mathbb{N}$ such that for all
$p\ge p_o$,
%
%
\begin{equation}
\label{ARpzeroes}
\inf_{|z|\le1+{1/p}} \Biggl| 1- \sum_{k=1}^p a_k(p)z^k \Biggr|
\ge\delta>0 .
\end{equation}
\end{lemma}

The uniform convergence of $A_p(z)$ toward $A(z)$ on the closed unit
disk immediately implies the following corollary to Lemma \ref{lemma3}.
\begin{corollary}
\label{corollary1}
Under the assumption of Lemma \ref{lemma3}, we have
%
%
\begin{equation}
\label{Azeroes}
A(z)= 1- \sum_{j=1}^{\infty} a_jz^j \neq0\qquad \forall|z|\le
1 .
\end{equation}
\end{corollary}

Lemma \ref{lemma3} and Corollary \ref{corollary1} now enable us to
invert the power series~$A(z)$ as well as the polynomial
$A_p(z)$. Let us denote
%
%
\begin{equation}
\label{ARinvert}
\Biggl( 1- \sum_{j=1}^{\infty} a_jz^j \Biggr) ^{-1} = 1+ \sum
_{j=1}^{\infty} \alpha_j z^j\qquad \forall|z| \le1
\end{equation}
and for all $p$ large enough (because of Lemma \ref{lemma3})
%
%
\begin{equation}
\label{ARpinvert}
\Biggl( 1- \sum_{j=1}^{p} a_j(p) z^j \Biggr) ^{-1} = 1+ \sum
_{j=1}^{\infty} \alpha_j (p) z^j\qquad \forall|z| \le1+\frac
{1}{p} .
\end{equation}
From (\ref{ARpinvert}), one immediately obtains that
%
%
\begin{equation}
| \alpha_j(p) | \le C\cdot\biggl( 1+\frac{1}{p}\biggr)
^{-j}\qquad \forall j \in\N.
\end{equation}
A further auxiliary result contains the transfer of the approximation
property of $a_j(p)$ for $a_k$ to the respective coefficients
$\alpha_k(p)$ and $\alpha_k$ of the inverted series. For such a
result, we make use of a weighted version of Wiener's lemma; cf.
\citet{Grochenig07}.
\begin{lemma}
\label{lemma4}
Under the assumptions of Lemma \ref{lemma3} and additionally $\sum
_{h=0}^{\infty} (1+h)^r |\gamma_X(h)|< \infty$ for some $r\ge0$
there exists a constant $C>0$ such that for all $p$ large enough
%
%
\begin{equation}
\label{alphap}
\sum_{j=1}^{\infty} (1+j)^r|\alpha_j(p)-\alpha_j| \le
C\cdot\sum_{j=p+1}^{\infty} (1+j)^r |a _j| \rightarrow_{p\to\infty}
0 .
\end{equation}
\end{lemma}

In a final step of this section, we now move on to estimators of the
coefficients $a_k(p)$. The easiest one might think of is to replace in
(\ref{ARplineareq}) the theoretical autocovariance function by its
sample version $\widehat\gamma_X(h)$.
Denote the resulting Yule--Walker estimators of $a_k(p)$ by $\widehat
a_k(p), k=1,\ldots, p$.
These Yule--Walker estimators are under the typically satisfied
assumption that $\widehat\gamma_X(0)>0$
uniquely determined and moreover fulfill (by the same arguments already
used) the desired property that
%
%
\begin{equation}
\label{zeroesAhat}
\widehat A _p(z)=1-\sum_{k=1}^p \widehat a_k(p) z^k \neq0\qquad
\forall|z|\le1 .
\end{equation}
Thus, we can also invert the polynomial $\widehat A_p(z)$ and we denote
%
%
\begin{equation}
\label{ARpdachinvert}
\Biggl( 1-\sum_{k=1}^p \widehat a_k(p) z^k\Biggr) ^{-1} = 1+\sum
_{k=1}^{\infty} \widehat\alpha_k(p) z^k ,\qquad |z|\le1 .
\end{equation}

We require that the estimators $(\widehat a_k(p):k=1,\ldots,p)$
converge---even at a~slow rate---to their theoretical counterparts,
namely:\vspace*{8pt}

(A1) $p(n)^2\cdot\sum_{1\le k\le p(n)} | \widehat
a_k(p(n)) - a_k(p(n)) | = {\mathcal O}_P(1)$, where $p(n)$
denotes a sequence of integers converging to infinity at a rate to be
specified.\vspace*{8pt}

Assumption (A1), for example, is met if a sufficient fast rate of
convergence for the empirical autocovariances toward their theoretical
counterparts can be guaranteed. The convergence\vadjust{\goodbreak} property of $\widehat
a_k(p)$ carries over to the corresponding coefficients
$\widehat\alpha_k(p)$ of the inverted polynomials [cf. (\ref
{ARpdachinvert})] as is specified in the following lemma.
\begin{lemma}
\label{lemma5}
Under the assumptions of Lemma \ref{lemma3} and \textup{(A1)}, we have
uniformly in $ k\in\mathbb{N}$
%
%
\begin{equation}
\label{alphapdach}
|\widehat\alpha_k(p(n))-\alpha_k(p(n))| \le
\biggl( 1 + \frac{1}{p(n)}\biggr) ^{-k} \frac{1}{p(n)^2} {\mathcal
O}_P(1) .
\end{equation}
\end{lemma}

\section{Validity of the AR-sieve bootstrap}
\label{bootfctlinstat}

\subsection{Functions of generalized means}
Consider a general class of estimators
%
%
\begin{equation}
\label{kunschstat}
T_n = f\Biggl( \frac{1}{n-m+1} \sum_{t=1}^{n-m+1} g( X_t,\ldots
,X_{t+m-1}) \Biggr) ,
\end{equation}
discussed in \citet{Kunsch89}, cf. Example 2.2; here $g\dvtx\mathbb
{R}^m \to\mathbb{R}^d$ and \mbox{$f\dvtx\mathbb{R}^d \to\mathbb{R}$}. For this
class of statistics, \citet{Buhlmann97} proved validity
of the AR-sieve bootstrap under the main assumption of an invertible
linear process with i.i.d. innovations for the underlying process
$(X_t)$; this means a process which admits an autoregressive
representation (\ref{eqarinfty}). The class of statistics
given in (\ref{kunschstat}) is quite rich and contains, for example,
versions of sample autocovariances, autocorrelations, partial autocorrelations,
Yule--Walker estimators and the standard sample mean as well.

The necessary smoothness assumptions on the functions $f$ and $g$ are
described below; these are identical
to the ones imposed by \citet{Buhlmann97}.\vspace*{8pt}

(A2) $f(y)$ has continuous partial derivatives for all $y$ in
a neighborhood of $\theta= E g(X_t,\ldots, X_{t+m-1})$ and the differentials
$\sum_{i=1}^m \partial f(y)/\partial x_i \vert_{x=\theta}$ do not vanish.
The function $g$ has continuous partial derivatives of order
$h$ $(h\ge1)$ that satisfy a Lipschitz condition.\vspace*{8pt}

We intend to investigate in this section what the autoregressive sieve
bootstrap really mimics if it is applied to statistics of the form
(\ref{kunschstat}) and the observations $X_1,\ldots,X_n$ do \textit{not}
stem from a linear AR($\infty$) process (\ref{eqarinfty}).
To be precise, we only assume that we
observe $X_1,\ldots,X_n$ from a process satisfying the following
assumption (A3).\vspace*{8pt}

(A3) $(X_t\dvtx t\in\mathbb{Z})$ is a zero mean, strictly
stationary and purely nondeterministic stochastic process. The
autocovariance function $\gamma_X$ satisfies\break
\mbox{$\sum_{h=0} ^{\infty} h^r \vert\gamma_X(h)\vert< \infty$} for
some $r\in\mathbb{N}$ specified in the respective results and
the spectral density $f_X$ is bounded and strictly positive.
Furthermore,\break $E(X_t^4)< \infty$.\vspace*{8pt}

Notice that the processes class described by (A3) is large enough
and includes several of the
commonly used linear and nonlinear time series models including
stationary and invertible autoregressive moving-average (ARMA) processes,
ARCH processes, GARCH\vadjust{\goodbreak} processes and so on. Summability of the
autocovariance function implies that the spectral density $f_X$ exists,
is bounded and continuous. We also added in (A3) the assumption
of finite fourth order moments of the time series. This assumption
seems to be unavoidable due to the autoregressive parameter estimation
involved in Step~1 of the AR-sieve bootstrap procedure and in regard
of assumption (A1).

From Section \ref{generalAR}, we know that if $(X_t)$ satisfies
assumption (A3) then it possesses an autoregressive
representation with an
uncorrelated white noise process $(\ve_t)$: cf.~(\ref
{woldar}). Because of the strict stationarity of $(X_t)$, we have that
the time series $(\ve_t)$
is strictly stationary as well and thus that the marginal distribution
$\mathcal{L}(\ve_t)$ of $ \ve_t$ does not depend on $t$.

Theorem \ref{bootlinstat} is the main result of this section. To state
it, we define the \textit{companion} autoregressive process
$\widetilde{\mathbf X}=({\widetilde X}_t\dvtx t \in\mathbb{Z})$ where $
\widetilde{X}_t$ is generated as follows:
%
%
\begin{equation}
\label{ARcompanion}
\widetilde X_t = \sum_{j=1}^{\infty} a_j \widetilde X_{t-j}
+\widetilde\ve_t,\qquad t\in\mathbb{Z} ;
\end{equation}
here $(\widetilde\ve_t)$ consists of i.i.d. random variables whose
marginal distribution of $ \widetilde\ve_t$ is identical
to that of $\ve_t $ from (\ref{woldar}), that is, $\mathcal
{L}(\widetilde\ve_t)=\mathcal{L}(\ve_t)$.
It is worth mentioning that all second order properties of $(\widetilde
X_t)$ and $(X_t)$, like autocovariance function and spectral density,
coincide while the probabilistic characteristics beyond
second order quantities
of both stationary processes are not necessarily the same.
Now, let $\widetilde T_n$ be the same statistic as $T_n$ defined in
(\ref{kunschstat}) but with
$X_t$ replaced by $\widetilde X_t$, that is,
%
%
\begin{equation}
\label{kunschstatcomp}
\widetilde{T}_n = f\Biggl( \frac{1}{n-m+1} \sum_{t=1}^{n-m+1} g(
\widetilde{X}_t,\ldots,\widetilde{X}_{t+m-1}) \Biggr).
\end{equation}

The main message of Theorem \ref{bootlinstat} is that the AR-sieve
bootstrap applied to data $X_1,\ldots,X_n$
in order to approximate the distribution of the statis\-tic~(\ref
{kunschstat}) will generally\vspace*{1pt} lead to an asymptotically
consistent estimation of the distribution of the statistic $\widetilde T_n$.
This implies that for the class of statis\-tics~(\ref{kunschstat}), the
AR-sieve bootstrap
will work \textit{if and only if} the limiting distributions of $ T_n$ and
of $\widetilde T_n$ are identical.
\begin{theorem}
\label{bootlinstat}
Assume assumptions (\textup{A1}), (\textup{A2}), (\textup{A3}) for
$r=1$ and the
moment condition
$E \ve_t^{2(h+2)}$ [cf. (\textup{A2}) for the definition of $h$ and
(\ref{woldar}) for the definition
of $\ve_t$], the condition $p(n)=o((n/\log n)^{1/4})$ on the
order of the approximating autoregression to the data and the following
two further assumptions:
\begin{longlist}[(A5)]
\item[(A4)] The empirical distribution function $F_n$ of the
random variables 
$ \ve_1,\ldots,\allowbreak \ve_n$ converges weakly to the distribution function
$F$ of $\mathcal{L}(\ve_1)$.
\item[(A5)] The empirical moments $1/n \sum_{t=1}^n\ve_t^r$
converge in probability to $E \ve_1^r$ for all $r\le2(h+2)$.
\end{longlist}
Then,
%
%
\begin{equation}
d_K \bigl( \mathcal{L}^\ast\bigl(\sqrt{n}\bigl( T_n^{\ast} - f(\theta
^{\ast})\bigr) \bigr) , \mathcal{L}\bigl( \sqrt{n}\bigl(\widetilde T_n -
f(\widetilde\theta) \bigr)\bigr) \bigr) =o_P(1)
\end{equation}
as $n\to\infty$.
Here $\theta^{\ast}= E g(X_t^{\ast},\ldots, X_{t+m-1}^{\ast})$,
$\widetilde\theta= E g(\widetilde X_t,\ldots, \widetilde
X_{t+m-1})$ and
$d_K$ denotes the Kolmogorov distance.
\end{theorem}

Some remarks are in order.
\begin{remark}
(i) Assumption (A4) does imply that we need some conditions
on the dependence structure of the random variables
$\ve_t$. For instance, a~standard mixing condition on $(\ve_t)$
suffices to ensure (A4): cf. \citet{PolitisRomanoWolf99},
Theorem 2.1.

\mbox{\hphantom{i}}(ii) Assumption (A5) on the empirical moments is fulfilled
if we ensure that sufficiently high empirical moments of the underlying
strictly stationary time series $X_t$ itself would converge in
probability to their theoretical counterparts.

(iii) As we already pointed out, Theorem \ref{bootlinstat}
states that the AR-sieve bootstrap mimics the behavior of the companion
autoregressive
process $(\widetilde X_t)$ as far as statistics of the form (\ref
{kunschstat}) are considered.
Of course if $(X_t)$ is a linear process with i.i.d. innovations and
if the corresponding moving average representation is invertible
leading to an infinite order autoregression (\ref{eqarinfty}) with
i.i.d. innovations, then the AR-sieve bootstrap works asymptotically as is
already known.
Moreover, for general process satisfying assumption~(A3), if we
are in the advantageous situation that the existing dependence
structure of the innovations $(\ve_t)$ appearing in the general
AR($\infty$) representation~(\ref{woldar}) does not show up in the
limiting distribution of~$ T_n$, then Theorem~\ref{bootlinstat} implies that the AR sieve
bootstrap works. We will illustrate this point by several examples
later on.
\end{remark}

Now we discuss relevant specializations of Theorem \ref{bootlinstat}.
Notice that the advantage of this theorem is that in order to check
validity of the AR-sieve bootstrap,
one only needs to check whether the
asymptotic distributions of $T_n=T_n(X_1,\ldots, X_n)$ based on the
observed time series is identical to the distribution
of the statistic $\widetilde T_n=T_n(\widetilde X_1,\ldots, \widetilde
X_n)$ based\vspace*{1pt} on fictitious observations $ \widetilde{X}_1,\widetilde
{X}_2, \ldots,
\widetilde{X}_n$ from the companion process $ \widetilde{\mathbf X}$.
If and only if this is the case, then the AR-sieve bootstrap works
asymptotically.
\begin{example}[(Sample mean)]
\label{mean}
Consider the case of the sample mean $ T_n=\overline{X}_n \equiv
n^{-1}\sum_{t=1}^nX_t$ and recall that under standard and
mild regularity conditions (e.g., mixing or weak dependence) we
typically obtain that
the sample mean of stationary time series satisfies $\sqrt{n}T_n
\Rightarrow N(0,\sum_{h=-\infty}^{\infty} \gamma_X(h))$ as $ n
\rightarrow\infty$ where
$ \Rightarrow$ stands for weak convergence.
Thus, the asymptotic distribution of the sample mean
depends only on the second order properties of the underlying process $
{\mathbf X}$ and since the companion
process $ \widetilde{\mathbf X}$ has the same second order properties
as $
{\mathbf X}$ we immediately get by Theorem
\ref{bootlinstat} that the AR-sieve bootstrap asymptotically works in
the case of the mean for general stationary time series for which
the spectral density is strictly positive. Even the strict stationarity
is not necessary in this case.
This is a~novel and significant extension of the results of \citet
{Buhlmann97}.
\end{example}
\begin{example}[(Sample autocovariances)]
\label{bootautocov} For $ 0 \leq h < n$, let $ T_n=\widehat{\gamma}(h)
\equiv n^{-1}\sum_{t=1}^{n-h}(X_t-\overline{X}_n)(X_{t+h}-\overline
{X}_n)$ be the sample autocovariance
at lag $ h$. Let us assume that
$ \sum_{h_1,h_2,h_3=- \infty}^{\infty}|\mathrm{cum} (X_t,X_{t+h_1}, X_{t+h_2},
X_{t+h_3} )| < \infty$ holds and denote by
\[
f_4(\lambda_1,\lambda_2,\lambda_3) = \sum_{h_1,h_2,h_3=-\infty}^{\infty
} \mathrm{cum} (X_t,X_{t+h_1}, X_{t+h_2}, X_{t+h_3} )e^{\imag \sum_{r=1}^{3}\lambda
_r h_r}
\]
the fourth order cumulant spectrum of ${\mathbf X}$. Under standard and
mild regularity conditions [see,
for instance, Dahlhaus (\citeyear{Dahlhaus85}), Theorem 2.1 and
Taniguchi and Kakizawa (\citeyear{TanKa00}), Chapter 6.1], it is known that
$ \sqrt{n}( T_n- \gamma(h)) \Rightarrow N(0,\tau^2)$ where
\begin{eqnarray*}
\tau^2 &=& 2\pi\int_{-\pi}^{\pi} 4\cos(\lambda h)^2 f^2(\lambda
)\,d\lambda\\
&&{}+ \iint4\cos(\lambda_1 h) \cos(\lambda_2 h)
f_4(-\lambda_1,\lambda_2,-\lambda_2)\,d\lambda_1 \,d\lambda_2.
\end{eqnarray*}
Notice that in contrast to the case of the sample mean, the limiting
distribution of the sample autocovariance depends also on the fourth
order moment structure of the underlying process $ {\mathbf X}$.
Now, to check validity of the AR-sieve bootstrap we have to derive the
asymptotic distribution of the sample autocovariances for the companion
autoregressive process $(\widetilde X_t)$.
This can be easily done, since the autoregressive polynomial in the
general autoregressive representation (\ref{woldar}) is always
invertible (cf.
Corollary \ref{corollary1}) from which we immediately get a one-sided
moving average representation with i.i.d. innovations
$(\widetilde\ve_t)$ for the companion
process $(\widetilde X_t)$. Furthermore, the fourth order cumulant
spectrum of $(\widetilde X_t)$ is given by
\[
\widetilde f_4(\lambda_1,\lambda_2,\lambda_3)=(2\pi)^{-3}\biggl( \frac
{E \ve_1^4}{(E \ve_1^2)^2}-3\biggr) \alpha(\lambda_1)\alpha(\lambda
_2)\alpha(\lambda_3)\alpha(-\lambda_1-\lambda_2-\lambda_3),
\]
where $ \alpha(\lambda)=\sum_{j=0}^{\infty} \alpha_j \exp\{-\imag j \lambda
\}$ and the coefficients $ (\alpha_k)$
are those appearing in (\ref{ARinvert}); see Section \ref{generalAR}.
Thus, for the sample autocovariance $ \widetilde{T}_n$ we get\vspace*{1pt} from
Brockwell and Davis [(\citeyear{BrockwellDavis91}), Proposition 7.3.1], that
$ \sqrt{n}( \widetilde T_n- \gamma(h)) \Rightarrow N(0,\widetilde{\tau
}^2)$ where
%
%
\begin{eqnarray}
\label{varianceofvariancetilde}
\widetilde{\tau}^2 & = & 2\pi\int_{-\pi}^{\pi} 4\cos(\lambda h)^2
f^2(\lambda)\,d\lambda\nonumber\\
&&{} + \iint4\cos(\lambda_1 h) \cos(\lambda_2 h)
\widetilde f_4(-\lambda_1,\lambda_2,-\lambda_2)\,d\lambda_1 \,d\lambda
_2 \nonumber\\[-8pt]\\[-8pt]
& = & \biggl( \frac{E \ve_1^4}{(E \ve_1^2)^2}-3\biggr)
(\gamma(h))^2\nonumber\\
&&{}+\sum_{k=-\infty}^{\infty} \bigl(\gamma(k)^2 + \gamma
(k+h) \gamma(k-h)\bigr).\nonumber
\end{eqnarray}
Since the variances $ \tau^2$ and $ \widetilde\tau^2$ of the
asymptotic distributions
of $ T_n$ and $\widetilde T_n$ do not coincide in general, we conclude
by Theorem \ref{bootlinstat} that the AR-sieve bootstrap fails for
sample autocovariances.
Notice that this failure is due to the fact that in general the
limiting distribution of sample autocovariances depends additionally on the
fourth order moment structure $f_4$ of the underlying process ${\mathbf
X}$, and this structure may substantially differ from that
of the companion process $\widetilde{\mathbf X}$.

Interestingly enough, even if the underlying process ${\mathbf X}$
is a linear time series, that is, satisfies (\ref{MAinfinity}),
the AR-sieve bootstrap may fail for the sample
autocovariances. To see why, note that
from the aforementioned proposition of \citet{BrockwellDavis91},
the asymptotic distribution of $T_n$ satisfies
$ \sqrt{n}( T_n- \gamma(h)) \Rightarrow N(0,\tau^2_L)$
where $\tau_L^2$ is given by
%
%
\begin{equation}
\label{varianceofvariance}
\tau_L^2= \biggl( \frac{E e_1^4}{(E e_1^2)^2}-3\biggr)
\gamma^2(h)+\sum_{k=-\infty}^{\infty} \bigl(\gamma^2(k) + \gamma
(k+h) \gamma(k-h)\bigr).
\end{equation}
%
Special attention is now due to the factor of
the first summand of (\ref{varianceofvariance}), which is the fourth
order cumulant of the i.i.d. process $(e_t)$.
Recall the asymptotic distribution of the sample autocovariances for
the companion autoregressive process $(\widetilde X_t)$
and especially its variance given in (\ref{varianceofvariancetilde}).
The two asymptotic variances given in (\ref{varianceofvariance}) and
(\ref{varianceofvariancetilde})
are in general not the same since the fourth order cumulant of the two
innovation processes $ (e_t)$ and
$ (\ve_t)$ are not necessarily the same.
We refer to Remark~\ref{remark1} for an example.
Of course, all appearing autocovariances are identical since we do not
change the second order properties of the process when
switching from $(X_t)$ to $(\widetilde X_t)$. The consequence is that
for sample autocovariances, the AR-sieve bootstrap
generally does not work even for linear processes of the type (\ref
{MAinfinity}) and this is true even if the process
is causal; see Remark \ref{remark1} and example (\ref{MA1}).
\begin{remark}
\label{rePoskitt}
If the innovations $ \{e_t\}$ in (\ref{MAinfinity}) are not necessarily
i.i.d. but form
a martingale difference sequence, then we are in the above described situation
more general than in the case of (\ref{MAinfinity}) with i.i.d.
innovations and thus the limiting distribution of sample
autocovariances and also of statistics of the type
(\ref{kunschstat}),
is not correctly mimicked by the
AR-sieve bootstrap. This contradicts Theorem 2 of \citet{Poskitt07}.
\end{remark}

We conclude this example by mentioning that in the case where the
i.i.d. innovations $e_t$ in (\ref{MAinfinity}) are normally distributed,
it follows that the random variables $\ve_t$ are
Gaussian as well, and that in both expressions (\ref
{varianceofvariance}) and~(\ref{varianceofvariancetilde})
the fourth order cumulants appearing as factors in the first summands vanish.
This means that for the special case of Gaussian time series fulfilling
assumption
(\ref{MAinfinity}) the AR-sieve bootstrap works.
\end{example}
%
%
\begin{example}[(Sample autocorrelations)] \label{bootautocorr}
Consider the estimator $ T_n= \widehat{\rho}(h)\equiv
\widehat{\gamma}(h)/\widehat{\gamma}(0)$ of the autocorrelation $ \rho_X(h)=
\gamma_X(h)/\gamma_X(0)$.
Due to the fact that for general processes satisfying assumption (A3)
the limiting distribution\vadjust{\goodbreak} of $T_n$ depends
also on the fourth order moment structure of the underlying process~$
{\mathbf X}$, cf. Theorem 3.1 of \citet{RomTho96},
which is not mimicked correctly by the companion process $ \widetilde
{\mathbf X}$, the AR-sieve bootstrap fails.
Fortunately the situation for autocorrelations is much better if we
switch to linear processes of type (\ref{MAinfinity}).
From Theorem 7.2.1 of \citet{BrockwellDavis91},
we obtain that $ \sqrt{n}(T_n-\rho_X(h)) \Rightarrow N(0,v^2)$ where
the asymptotic variance is given by Bartlett's formula:
\[
v^2=\sum_{k\in\Z}\bigl\{\bigl(1+2\rho^{2}_X(h)\bigr)\rho^2_X(k)+\rho_X(k-h)\rho
_X(k+h) - 4\rho_X(h)\rho_X(k)\rho_X(k+h) \bigr\}.
\]
As can be seen, in this case the asymptotic variance depends only on the
autocorrelation function $ \rho_X(\cdot)$ [or equivalently on the
standardized spectral density $ \gamma_X^{-1}(0)f_X(\cdot)$] of the
underlying process $ {\mathbf X}$.
This means that the first summand in (\ref{varianceofvariance}) which
refers to the fourth order cumulant of the i.i.d. innovation process
$e_{t}$ in (\ref{MAinfinity})
does not show up in the limiting distribution of sample autocorrelations
and this in turn leads to the fact that the asymptotic distribution of
sample autocorrelations based on observations stemming from the process
${\mathbf X}$ is identical to that of the sample autocorrelation based on
observations stemming from the companion autoregressive process
$\widetilde{\mathbf X}$. This is true
since the
companion autoregressive process $\widetilde{\mathbf X}$ shares all second
order properties of the process ${\mathbf X}$. Hence,
the AR-sieve bootstrap works for the autocorrelations given data from
a linear process (\ref{MAinfinity}). We stress here
the fact that this result is
true regardless whether the representation (\ref{MAinfinity}) allows
for an autoregressive inversion or not and holds even though the
probabilistic properties of the underlying process $(X_t)$ and of the
autoregressive companion process $(\widetilde X_t)$ beyond second order
properties are not the same.
\end{example}
\begin{remark}
\label{othestat}
Similar to the autocorrelation case,
the validity of the AR-sieve bootstrap in the linear class (\ref{MAinfinity})
is shared by many different statistics whose large-sample distribution
depends only on the (first and) second order moment structure of the
underlying process.
Examples include the partial autocorrelations or Yule--Walker estimators
of autoregressive coefficients.
\end{remark}
\begin{remark}
\label{runningmean}
Consider statistics of the type (\ref{kunschstat}) in the easiest case
where $d=1$, that is,
%
%
\begin{equation}
\frac{1}{n-m+1} \sum_{t=1}^{n-m+1} g(X_t,\ldots,X_{t+m-1})
\end{equation}
for a function $g\dvtx\mathbb{R}^m\to\mathbb{R}$. Here, the practitioner
may approach this
as the sample mean of observations $Y_1, \ldots, Y_{n-m+1}$ where
$(Y_t=g(X_t,\ldots,X_{t+m-1})\dvtx\allowbreak t\in\mathbb{Z})$. Notice that strict
stationarity as well as mixing properties easily carries over from
$(X_t)$ to $(Y_t)$. Thus, we may apply the AR-sieve bootstrap to the
sample mean of $Y_t$, which works under quite general
assumptions; cf. Remark~\ref{mean} and Theorem~\ref{bootlinstat}. The
only\vadjust{\goodbreak} crucial assumption for establishing
asymptotic consistency of the AR-sieve bootstrap is that we need the
property that the spectral density $f_Y$ of the transformed time series $(Y_t)$
is strictly positive and continuous. Although this is not a very
restrictive condition, it may be difficult to check since
there seems to be no general result describing the behavior of spectral
densities
under nonlinear transformations. 
\end{remark}

\subsection{Integrated periodograms}
The considerations of Examples \ref{bootautocov} and~\ref{bootautocorr}
can be transferred to integrated periodogram estimators for these
quantities which lead us to the second large class of statistics that
we will discuss.
Denote, based on observations $X_1,\ldots,X_n$, the periodogram
$I_n(\lambda)$ defined~by
%
%
\begin{equation}
\label{periodogram}
I_n(\lambda)= \frac{1}{2\pi n}\Biggl\vert\sum_{t=1}^n X_t e^{-\imag
\lambda t} \Biggr\vert^2 ,\qquad \lambda\in[0,\pi],
\end{equation}
and consider a general class of integrated periodogram estimators defined~by
%
%
\begin{equation}
\label{integratedperiodogram}
M( I_n,\varphi) = \int_0^{\pi} \varphi(\lambda
)I_n(\lambda) \,d\lambda,
\end{equation}
where $\varphi$ denotes an appropriately defined function on $[0, \pi
]$. Under the main assumption that the underlying process $ {\mathbf X}$
has
the representation (\ref{MAinfinity}), \citet{Dahlhaus85}
investigated asymptotic properties of (\ref{integratedperiodogram}) and
obtained the asymptotic distribution of $\sqrt{n}(M(I_n,\varphi
)-M(f_X,\varphi))$. In particular,
it has been shown that $\sqrt{n}(M(I_n,\varphi)-M(f_X,\varphi))$
converges, as $ n \rightarrow\infty$,
to a~Gaussian distribution with zero mean and variance given by
%
%
\begin{equation}
\label{varianceintperiodo}
\bigl(E e_1^4/(E e_1^2)^2 -3\bigr)\biggl( \int_0^{\pi} \varphi(\lambda
)f_X(\lambda) \,d\lambda\biggr) ^2 + 2\pi\int_0^{\pi} \varphi
^2(\lambda)
f_X^2(\lambda) \,d\lambda.
\end{equation}
%
Notice that substituting $\varphi(\lambda)$ by $2\cos(\lambda h)
, h\in\mathbb{N}_0$, implies that $M(I_n,\varphi)$
equals the sample autocovariance of the observations $X_1,\ldots,X_n$
at lag $h$ and (\ref{varianceintperiodo}) would then exactly
turn to be the asymptotic variance given in~(\ref{varianceofvariance}).


As we will see in Theorem \ref{intperiodoCLT}, the situation for
integrated periodo\-grams~(\ref{integratedperiodogram})~is rather similar
to that
of empirical autocovariances which are of course special cases of
integrated periodograms. Thus,
we only discuss briefly this rather relevant class of statistics. As
Theorem \ref{intperiodoCLT} shows,
we obtain for this class that the AR-sieve bootstrap
asymptotically mimics the behavior of $\sqrt{n}(M(\widetilde I_n,
\varphi)-M(f_X,\varphi))$, where $\widetilde I_n$ is defined\vspace*{1pt} as
$I_n $ with~$X_t$ replaced by the companion autoregressive time series
$\widetilde X_t$, that is,
%
%
\begin{equation}
\label{periodogramtilde}
\widetilde{I}_n(\lambda)= \frac{1}{2\pi n}\Biggl\vert\sum_{t=1}^n
\widetilde{X}_t e^{-\imag\lambda t} \Biggr\vert^2 ,\qquad \lambda\in
[0,\pi] .
\end{equation}
%
%
\begin{theorem}
\label{intperiodoCLT}
Assume \textup{(A1)} and \textup{(A3)} with $r=1$ and assume that for all
$M\in\mathbb{N}$
%
%
\begin{equation}
\label{gammaCLT}
{\mathcal L}^\ast\bigl( \sqrt{n}\bigl( \widehat\gamma_{X^{\ast}}(h)-
E^{\ast} \widehat\gamma_{X^{\ast}(h)} \bigr) \dvtx
h=0,1,\ldots,M \bigr)
\Rightarrow N( 0, V_M )
\end{equation}
(in probability), where $\widehat\gamma_{X^{\ast}}(h)$ denote the
empirical autocovariances of the bootstrap observations $ X_1^\ast,
X_2^\ast, \ldots, X_n^\ast$ and where
%
%
\begin{eqnarray}
\label{asymptotVAR}
V_M&=& \Biggl[ \biggl( \frac{E \ve_1^4}{(E \ve_1^2)^2} -3 \biggr)
\gamma(i)\gamma(j) \nonumber\\[-8pt]\\[-8pt]
&&\hspace*{4.7pt}{}+\sum_{k=-\infty}^{\infty} \bigl(\gamma(k)
\gamma(k-i+j)+\gamma(k+j)\gamma(k-i)\bigr) \Biggr] _{i,j=0}^M .\nonumber
\end{eqnarray}
Then we obtain for all $\varphi$ bounded and with bounded variation
that (in probability)
%
%
\begin{eqnarray}
\label{bootintperiodo}
&&d_K \bigl( \mathcal{L}^\ast\bigl( \sqrt{n}\bigl(M(I_n^{\ast},\varphi
)-M(\widehat f_{\mathrm{AR}},\varphi)\bigr)\bigr)
,\nonumber\\[-8pt]\\[-8pt]
&&\qquad\hspace*{10pt}
\mathcal{L} \bigl( \sqrt{n}\bigl( M(\widetilde I_n,\varphi)-M(f_X,\varphi)
\bigr) \bigr)\bigr)
\to0 .\nonumber
\end{eqnarray}
Here $\widehat f_{\mathrm{AR}}(\lambda)=\frac{\widehat\sigma
^2(p(n))}{2\pi} \vert1-\sum_{j=1}^{p(n)}
\widehat a_j(p(n))e^{-\imag j\lambda} \vert^{-2} , \lambda
\in[0,\pi]$
and $\widehat\sigma(p(n))^2=E \widehat\ve_1(p(n))^2$, cf.
Step 2 in the definition of the AR-sieve bootstrap procedure.

Moreover, the limiting Gaussian distribution of $ \sqrt{n}(M(I_n^{\ast
},\varphi)-M(\widehat f_{\mathrm{AR}},\varphi)$ possesses the following
variance:
%
%
\begin{equation}
\label{varperiodo}\qquad
\bigl(E \ve_1^4/(E \ve_1^2)^2 -3\bigr)\biggl( \int_0^{\pi} \varphi(\lambda
)f_X(\lambda) \,d\lambda\biggr) ^2 + 2\pi\int_0^{\pi} \varphi
^2(\lambda)
f_X^2(\lambda) \,d\lambda.
\end{equation}
%
\end{theorem}
\begin{remark}
\label{remarkintperiodoCLT}
(i) Assumptions under which (\ref{gammaCLT}) is fulfilled are
given in Theorem \ref{bootlinstat} since sample autocovariances
belong to the class (\ref{kunschstat}).

(ii) Theorem \ref{intperiodoCLT} implies that if $\int_0^{\pi
}\varphi(\lambda) f(\lambda) \,d\lambda= 0$ then the AR-sieve
bootstrap asymptotically
works
for the integrated periodogram statistics $M(I_n,\varphi)$ for time
series fulfilling (\ref{MAinfinity}). This follows immediately
by a comparison of the asymptotic variances (\ref{varianceintperiodo})
and (\ref{varperiodo}). Clearly, the
same result holds true if the underlying time series is normally
distributed since in this case
both innovation processes, $ (\ve_t)$ and $ (e_t)$, are Gaussian and
therefore the fourth order cumulants vanish.
In all other cases the AR-sieve bootstrap does not work in general,
since the fourth order cumulant $E \ve_1^4/(E
\ve_1^2)^2 -3$ does not necessarily coincide with $E e_1^4/(E
e_1^2)^2 -3$; see Example \ref{bootautocov}.
\end{remark}

Relevant statistics for which we can take advantage of the condition
$\int_0^{\pi}\varphi(\lambda)\times f(\lambda) \,d\lambda= 0$ are the
so-called ratio statistics which are defined by
%
%
\begin{equation}
\label{ratio}
R(I_n,\varphi) = \frac{M(I_n,\varphi)}{\int_0^{\pi} I_n(\lambda
) \,d\lambda} .
\end{equation}
For this class of statistics, \citet{DahlhausJanas96} showed that
under the same assumptions of a linear process of type
(\ref{MAinfinity}) one obtains that
$\sqrt{n}(R(I_n,\varphi)-R(f_X,\varphi))$ has a Gaussian limiting
distribution with mean zero and variance given by
%
%
\begin{eqnarray}
\label{varianceratio}
\frac{2\pi\int_0^{\pi} \psi^2(\lambda) f_X^2(\lambda) \,d\lambda
}{( \int f_X(\lambda) \,d\lambda) ^4}\hspace*{50pt}\nonumber\\[-8pt]\\[-8pt]
&&\eqntext{\mbox{where } \displaystyle \psi(\lambda) = \varphi(\lambda) \int f_X(\lambda)
\,d\lambda
-\int\varphi(\lambda) f_X(\lambda) \,d\lambda.}
\end{eqnarray}
Thus, exactly as in the case of sample autocorrelations the fourth
order cumulant term [cf. (\ref{varianceintperiodo})] of the i.i.d.
innovation process disappears and therefore again the following
corollary to Theorem \ref{intperiodoCLT} is true. This corollary
states that the AR-sieve bootstrap works for ratio statistics under the
quite general assumption
that the underlying process is a linear time series
(\ref{MAinfinity}) with i.i.d. innovations and a strictly positive
spectral density which under model (\ref{MAinfinity}) is always continuous.
\begin{corollary}
\label{intperiodoCLTcor}
Under the assumptions of Theorem \ref{intperiodoCLT}, we have that (in
probability)
%
%
\begin{equation}
\label{bootintperiodoratio}
d_K \bigl( \mathcal{L}^\ast\bigl( \sqrt{n}\bigl(R(I_n^{\ast},\varphi
)-R(\widehat f_{\mathrm{AR}},\varphi)\bigr)\bigr) ,
\mathcal{L} \bigl( \sqrt{n}\bigl( R(I_n,\varphi)-R(f,\varphi) \bigr)
\bigr)\bigr)
\rightarrow0 .\hspace*{-40pt}
\end{equation}
Moreover, the limiting Gaussian distribution of $ \sqrt{n}(R(I_n^{\ast
},\varphi)-R(\widehat f_{\mathrm{AR}},\varphi))$ possesses the variance
given in (\ref{varianceratio}).
\end{corollary}

Another class of integrated periodogram estimators is that of
nonparametric estimators of the spectral density $f_X$
which are obtained from (\ref{integratedperiodogram}) if we allow for
the function $ \varphi$ to depend on $n$. In particular, let $ \varphi
_n(\lambda)=K_h(\omega-\lambda)$
for some $\omega\in[0,\pi]$ where
$ h=h(n)$ is a sequence of positive numbers (bandwidths) approaching
zero as $n \rightarrow\infty$, $ K_h(\cdot)=h^{-1}K(\cdot/h)$ and $K$
is a kernel function satisfying
the following assumption:\vspace*{8pt}

(A6) $K$ is a nonnegative kernel function with compact support
$[-\pi,\pi]$. The Fourier transform $k$ of $K$ is assumed to be a
symmetric, continuous and
bounded function satisfying $k(0)=2\pi$ and $ \int_{-\infty}^{\infty
}k^2(u)\,du < \infty$.\vspace*{8pt}

Denote by $ f_{n,X} $ be the resulting integrated periodogram
estimator, that is,
%
%
\begin{equation} \label{specdenest}
f_{n,X}(\omega)= \int_{-\pi}^{\pi} K_h(\omega-\lambda) I_n(\lambda)
\,d\lambda.
\end{equation}
Notice that the asymptotic properties of the estimator (\ref
{specdenest}) of the spectral density are identical to those of its
discretized version
$ \widehat{f}_{n,X}(\omega)= \break (nh)^{-1} \sum_{j}K_h(\omega-\lambda_j)
I_n(\lambda_j) $, where $\lambda_j=2\pi j/n$ are the Fourier
frequencies, as well as of
so-called lag-window estimators; cf. \citet{Priestley81}.

Now, let $ f^\ast_{n,X}$ be the same estimator as (\ref{specdenest})
based on the AR-sieve bootstrap periodogram
$ I^\ast_n(\lambda)=(2\pi n)^{-1}|\sum_{t=1}^{n}X_t^\ast\exp\{\imag\lambda
t\}|^2 $. We then obtain the following theorem.
\begin{theorem}
\label{intperiodoCLTspec} Under the assumptions of Theorem \ref
{intperiodoCLT} with $r=2$ and assumption \textup{(A6)}, we have that (in probability)
%
%
\begin{equation}
\label{dKbootspec}
d_K \bigl( \mathcal{L}^\ast\bigl( \sqrt{nh}\bigl(f^{\ast}_{n,X}(\lambda)-
\widehat f_{\mathrm{AR}}(\lambda) \bigr)\bigr) ,
\mathcal{L} \bigl( \sqrt{nh}\bigl( f_{n,X}(\lambda)- f_{X}(\lambda) \bigr)
\bigr) \bigr) \rightarrow0 .\hspace*{-28pt}
\end{equation}
Moreover, conditionally on $ X_1, X_2, \ldots, X_n$,
%
%
\begin{eqnarray}
\label{biassd}
&&\sqrt{nh}E^\ast\bigl(f^{\ast}_{n,X}(\lambda)- \widehat f_{\mathrm
{AR}}(\lambda)
\bigr)\nonumber\\[-8pt]\\[-8pt]
&&\qquad \rightarrow\cases{
0, &\quad if $n^{-1/5}h \to0$, \vspace*{2pt}\cr
\displaystyle\frac{1}{4\pi} f_X^{\prime\prime}(\lambda)
\int u^2 K(u)\,du, &\quad if
$n^{-1/5}h \to1$,}\nonumber
\end{eqnarray}
where $ f_X^{\prime\prime}$ denotes the second derivative of $f_X$ and
%
%
\begin{equation}
\label{varsd}
nh \operatorname{Var} (f^{\ast}_{n,X}(\lambda)) \rightarrow( 1+\delta
_{0,\pi
}) f^2_X(\lambda) (2\pi)^{-1}\int K^{2}(u)\,du,
\end{equation}
where $ \delta_{0,\pi}=1 $ if $ \lambda=0 $ or $ \pi$ and $ \delta
_{0,\pi}=0$ otherwise.
\end{theorem}

Recall that under Assumption (A3) it has been shown under
different regularity conditions [see \citet{ShaoWu07}] that $ \sqrt
{nh}( f_{n,X}(\lambda)- f_{X}(\lambda) )$
converges to a Gaussian distribution with mean and variance given by
the expression on the right-hand side of (\ref{biassd}) and (\ref
{varsd}), respectively.
Thus, the above theorem implies that for spectral density estimators
like (\ref{specdenest}), the AR-sieve bootstrap asymptotically is valid
for a very broad class of stationary time series that goes far beyond
the linear processes class (\ref{MAinfinity}).

Corollary \ref{intperiodoCLTcor} and Theorem \ref{intperiodoCLTspec}
highlight an interesting relation between frequency domain bootstrap
procedures, like for instance those proposed by
\citet{FrankeHardle9299} and
\citet{DahlhausJanas96} and the AR-sieve bootstrap. Notice that
the basic assumptions imposed on the underlying process ${\mathbf X}$ for
such a frequency domain
bootstrap procedure to be valid are that the underlying process
satisfies (\ref{MAinfinity}) with a strictly positive spectral density $f_X$.
Furthermore, validity of such a frequency domain procedure has
been established only for those statistics
for which their limiting distribution does not depend on the fourth
order moment structure of the innovation process $e_{t}$ in (\ref{MAinfinity}).
Thus, such a frequency domain bootstrap essentially works for
statistics like ratio statistics \label{ratio} or nonparametric
estimators of the spectral density
like (\ref{specdenest}).
The results of this section, that is, Theorem~\ref{intperiodoCLT}, Corollary
\ref{intperiodoCLTcor} and Theorem \ref{intperiodoCLTspec}, imply that
if the underlying stationary process satisfies (\ref{MAinfinity}) and
if the spectral density is strictly positive then
the AR-sieve bootstrap works in the same cases in which the
frequency domain bootstrap procedures work.

\section{Conclusions}
\label{conclusions}

In this paper, we have investigated the range of validity of the
AR-sieve bootstrap. Based on a quite general Wold-type autoregressive
representation, we provided a simple and effective tool for verifying
whether or not the AR-sieve bootstrap asymptotically works. The central
question is to what extent the complex dependence structure of the
underlying stochastic process shows up in the (asymptotic) distribution
of the relevant statistical quantities. If the asymptotic behavior of
the statistic of interest based on our
data series 
is identical to that of the same statistic based on data generated from
the companion autoregressive process,
then the AR-sieve bootstrap leads to asymptotically correct results.

The family of estimators that have been considered ranges from simple
arithmetic means
and sample autocorrelations to quite general sample means of functions
of the observations as well as spectral density estimators and
integrated periodograms. Our concrete findings concerning validity of
the AR-sieve bootstrap are different for different statistics.
Generally speaking, if
the asymptotic distribution of a relevant statistic is determined
solely by the first and second order
moment structure, then the AR-sieve bootstrap is expected to work.
Thus, validity of the AR-sieve bootstrap does not require that the
underlying stationary process obeys a linear AR($\infty$)
representation (with i.i.d. or martingale difference errors) as was
previously thought. Indeed,
for many statistics of interest, the range of the validity of the
AR-sieve bootstrap goes far
beyond this subclass of linear processes. In contrast, we point out the
possibility that the AR-sieve bootstrap
may fail even though the data series \textit{is} linear; a prominent example
is the sample
autocovariance in the case of the data arising from a noncausal AR($p$)
or a~noninvertible~MA($q$) model.

Finally, our results bear out an interesting analogy between frequency
domain bootstrap methods and
the AR-sieve method. In the past, both of these methodologies have been
thought to work only in the linear time series setting.
Nevertheless, we have just shown the validity of
the AR-sieve bootstrap for many statistics of interest without the
assumption of
linearity, for example, under the general assumption (A3) and some
extra conditions.
In recent literature, some examples have been found where the frequency
domain bootstrap also works without the assumption of a linear process;
see, for example,
the case of spectral density estimators studied by
\citet{ShaoWu07}. By analogy to the AR-sieve results of the
paper, it can be conjectured that frequency
domain bootstrap methods might also be valid without the linearity assumption
as long as the statistic in question has a~large-sample
distribution depending only on first and second order moment properties;
cf. \citet{KirchPol2011} for some results in that
direction.


\begin{appendix}\label{proofs}
\section*{Appendix: Auxiliary results and proofs}

\vspace*{-6pt}

\begin{pf*}{Proof of Lemma \ref{lemma2}}
From Baxter [(\citeyear{Baxter62}), Theorem 2.2] in a slightly more general
version given in Baxter [(\citeyear{Baxter63}), Theorem 1.1], we obtain for
arbitrary submultiplicative weight or norm functions $\nu(k)\ge1$,
that is, $\nu(n) \le\nu(m)\cdot\nu(n-m)$ for all $n,m$, that the
following bound holds true
for all $p\in\mathbb{N}$ and a constant $C>0$
%
%
\begin{equation}
\label{Bax1}
\sum_{k=0}^{p} \nu(k) \biggl\vert\frac{a_k(p)}{\sigma^2(p)} -
\frac{a_k}{\sigma^2_{\varepsilon}} \biggr\vert
\le C\cdot\sum_{k=p+1}^{\infty} \nu(k) \biggl\vert\frac
{a_k}{\sigma_{\varepsilon}^2} \biggr\vert.
\end{equation}
Here, $\sigma^2(p) = E ( X_t-\sum_{k=1}^p a_k(p) X_{t-k}
) ^2 \le\sigma^2 (0)$ for all $p\in\mathbb{N}$
and $\sigma^2(p) \to\sigma_{\varepsilon}^2$ [cf. (\ref
{sigmasquare})] as $p\to\infty$.

Since $\nu(k)=(1+k)^r$ is submultiplicative for all $r\ge0$, cf.
Gr{\"o}chenig [(\citeyear{Grochenig07}), Lemma 2.1], we obtain from (\ref{Bax1})
%
%
\begin{eqnarray}
\label{Bax2}
&& \sum_{k=0}^p
(1+k)^r \vert a_k(p)-a_k\vert\nonumber\\
&&\qquad\le \sum_{k=0}^{p} (1+k)^r \biggl\vert\frac{a_k(p)}{\sigma
^2(p)} - \frac{a_k}{\sigma^2_{\varepsilon}} \biggr\vert
\cdot\sigma^2(p) \nonumber\\[-8pt]\\[-8pt]
&&\qquad\quad{}+ \sum_{k=0}^{p} (1+k)^r \vert a_k \vert
\biggl\vert\frac{1}{\sigma^2_{\varepsilon}} -
\frac{1}{\sigma^2(p)} \biggr\vert\cdot\sigma^2(p) \nonumber\\
&&\qquad\le \frac{C \sigma^2(0)}{\sigma^2_{\varepsilon}} \Biggl( 1+ \sum
_{k=0}^{\infty} (1+k)^r \vert a_k \vert\Biggr)
\cdot\sum_{k=p+1}^{\infty} (1+k)^r \vert a_k \vert, \nonumber
\end{eqnarray}
which is the assertion of Lemma \ref{lemma2}. To see the last bound,
observe that because of $a_0(p)=a_0=1$
we can bound $ \vert\frac{1}{\sigma^2(p)} - \frac{1}{\sigma
^2_{\varepsilon}} \vert$ by the right-hand side of (\ref
{Bax1}) as well.
\end{pf*}
\begin{pf*}{Proof of Lemma \ref{lemma3}} As mentioned just in front of
the statement of Lem\-ma~\ref{lemma3},
we have $A_p(z)=1-\sum_{k=1}^p a_k(p)z^k$ for all $|z|\le1$. Now
assume that
(\ref{ARpzeroes}) is false. Then there exists a sequence $\{ p(k)\dvtx
k\in
\mathbb{N}\} \subset\mathbb{N}$, $p(k)\to\infty$ and
a sequence $\{z_k\dvtx k\in\mathbb{N}\}$ of complex numbers with
$|z_k|\le
1+1/p(k)$ such that
%
%
\begin{equation}
\label{Apk}
A_{p(k)}(z_k) \to_ {k\to\infty} 0 .
\end{equation}
Let us further assume that we can find a subsequence of $\{z_k\dvtx k\in
\mathbb{N}\}$ which completely stays within the closed unit disk.
Without loss of generality, assume that $\{z_k\}$ itself has this
property. Since we have \mbox{$A_p(z)\neq0, \forall|z|\le1$} and because
$A_p$ is holomorphic, the minimum principle of holomorphic
functions\vadjust{\goodbreak}
leads to
%
%
\begin{equation}
| A_p(z) | \ge\min_{|z|=1} | A_p(z) | \qquad \forall
|z|\le1.
\end{equation}
The set $\{ z\in\mathbb{C}| |z|=1\}$ is compact and $|A_p|$ is
continuous, thus there exists a $z_p^{\ast}$ with $|z_p^{\ast}|=1$ and
$|A_p(z_p^{\ast})| = \min_{|z|=1} |A_p(z)|$.

Without loss of generality, assume that $z_p^{\ast}$ converges to a
complex number~$z^{\ast}$ with $|z^{\ast}|=1$. From the above, we have
%
%
\begin{equation}
\label{Apconv}
\bigl| A_{p(k)}\bigl(z_{p(k)}^{\ast}\bigr) \bigr| \le\bigl| A_{p(k)}(z_{k})
\bigr| \rightarrow_{k\to\infty} 0 .
\end{equation}
Writing
%
%
\begin{equation}\quad
A(z^{\ast})=A(z^{\ast})-A\bigl(z_{p(k)}^{\ast}\bigr) + A\bigl(z_{p(k)}^{\ast
}\bigr)-A_{p(k)}\bigl(z_{p(k)}^{\ast}\bigr) + A_{p(k)}\bigl(z_{p(k)}^{\ast}\bigr)
\end{equation}
and having in mind that $A_p(z)$ converges to $A(z)$ uniformly on the
closed unit disk because of Lemma \ref{lemma2}
and regarding (\ref{Apconv}) as well as the continuity of $A(z)$ we
finally obtain $A(z^{\ast})=0$ which is a contradiction to $A(z) \neq
0$ for all $\vert z \vert=1$ [cf. below (\ref{Azeroesfirst})].

Since we cannot find a subsequence of $z_k$, completely staying in the
unit disk it exists a subsequence
$(z_{k^{\prime}})$ that completely stays in the region
$1<|z|\le1+1/p(k^{\prime})$. Again assume without loss of generality
that $k^{\prime}=k$ and that $z_k$ converges to some
$z_o$ which necessarily must fulfill $|z_o|=1$.

We will show that $A(z_o)=0$ holds, which again is a contradiction to
$A(z) \neq0$ for all $\vert z \vert=1$.
To this end, let us write $A(z_o)$ in the following way:
%
%
\begin{eqnarray}
A(z_o)&=&A_{p(k)}(z_k) + \sum_{j=1}^{p(k)} \bigl( a_j(p(k))-a_j\bigr)
z_k^j\nonumber\\[-8pt]\\[-8pt]
&&{}+ \sum_{j=1}^{p(k)} a_j(z_k^j-z_o^j)
- \sum_{j=p(k)+1}^{\infty} a_j z_o^j .\nonumber
\end{eqnarray}
The first summand on the right-hand side converges to zero by (\ref
{Apk}) and the last summand is bounded through
$\sum_{j=p(k)+1}^{\infty} |a_j| \to0$ as $k\to\infty$. The second
summand in turn is bounded by
%
%
\begin{equation}
\sup_{|z|\le1+{1/p}} \Biggl| \sum_{j=1}^p \bigl(a_j(p)-a_j\bigr)z^j
\Biggr| \le
\sum_{j=1}^p |a_j(p)-a_j| \sup_{|z|\le1+{1/p}} |z|^p
\rightarrow_{p\to\infty} 0.\hspace*{-28pt}
\end{equation}
For the third and last summand, which reads
%
%
\begin{equation}
\sum_{j=1}^{\infty} a_j (z_k^j-z_o^j) 1\{ j\le p(k)\} ,
\end{equation}
one obtains by dominated convergence [recall that $|z_k|^j$ for $j\le
p(k)$ is bounded by $3$ and that $z_k\to z_o$] also convergence to zero.

This concludes the proof of Lemma \ref{lemma3}.\vadjust{\goodbreak}
\end{pf*}
\begin{pf*}{Proof of Lemma \ref{lemma4}}
Under the assumptions the autoregressive coefficients, $a_k$ have the
following property:
%
%
\begin{eqnarray}
\mathbf{a} &:=& (1,-a_1,-a_2, \ldots) \in\ell_1^v \nonumber\\[-8pt]\\[-8pt]
&:=& \Biggl\{ (z_j\dvtx j\in
\mathbb{N}_0)\subset\mathbb{C} \Big| \sum_{j=0}^{\infty}
(1+j)^r |z_j| < \infty\Biggr\} .\nonumber
\end{eqnarray}
Because of $1-\sum_{j=1}^{\infty}a_jz^j\neq0$ $\forall|z|\le1$ (cf.
Corollary \ref{corollary1}) we have from Gr{\"o}chenig
[(\citeyear{Grochenig07}), Theorem 6.2] that a multiplicative inverse
${\mathbf a}^{-1} \in\ell_1^v$ exists. For this result, observe that
our weight function $\nu(k)=(1+k)^r$ satisfies the so-called
Gelfand--Raikov--Shilov (GRS) condition
%
%
\begin{equation}
\label{GRS}
\nu(nk)^{1/n} \to1 \qquad\mbox{as } n \to\infty;
\end{equation}
cf. \citet{Grochenig07}, Lemma 2.1.

Since multiplication here is the usual convolution of sequences, we
have that
${\mathbf a}^{-1} =(1, \alpha_1,\alpha_2, \ldots)$, where the
coefficients $\alpha_k$ coincide with the power series coefficients of
$(1-\sum_{k=1}^{\infty} a_kz^k)^{-1}$. The assertion ${\mathbf a}^{-1}
\in\ell_1^v$ then just means that we have for the coefficients
$\alpha_k$ in
(\ref{ARpinvert})
%
%
\begin{equation}
\label{alpha}
\sum_{j=0}^{\infty} (1+j)^r |\alpha_j| < \infty.
\end{equation}
Exactly along the same lines, we obtain [cf. (\ref{ARpinvert})]
%
%
\begin{eqnarray}
{\mathbf a}(p)^{-1}&=&(1,-a_1(p),-a_2(p),\ldots, -a_p(p),0,\ldots)
^{-1}\nonumber\\[-8pt]\\[-8pt]
& = &(1,\alpha_1(p),\alpha_2(p),\ldots) \in\ell_1^v . \nonumber
\end{eqnarray}
We have
\begin{eqnarray*}
&& \sum_{k=0}^{\infty} (1+k)^r |\alpha_k(p)-\alpha_k| \\
&&\qquad= | {\mathbf a(\mathbf p)}^{-1} - {\mathbf a}^{-1} | _{\ell_1^v}
\\
&&\qquad= \bigl| {\mathbf a(\mathbf p)}^{-1}\bigl({\mathbf a} - {\mathbf a(\mathbf
p)}\bigr){\mathbf a}^{-1} \bigr|
_{\ell_1^v} \\
&&\qquad= | {\mathbf a(\mathbf p)}^{-1}-{\mathbf a}^{-1} | _{\ell_1^v}
| {\mathbf a}-{\mathbf a(\mathbf p)} | _{\ell_1^v} + | {\mathbf a}^{-1}
| _{\ell_1^v}
+ | {\mathbf a}^{-1} | _{\ell_1^v} | {\mathbf a}-{\mathbf a(\mathbf p)}
| _{\ell_1^v}.
\end{eqnarray*}
Simple algebra finally leads to
\[
\sum_{k=0}^{\infty} (1+k)^r |\alpha_k(p)-\alpha_k|
= | {\mathbf a(\mathbf p)}^{-1} - {\mathbf a}^{-1} | _{\ell_1^v}
\le\frac{ | {\mathbf a}^{-1} | _{\ell_1^v}^2 | {\mathbf a}-{\mathbf
a(\mathbf p)} | _{\ell_1^v}}
{1- | {\mathbf a}^{-1} | _{\ell_1^v} | {\mathbf a}-{\mathbf a(\mathbf p)}
| _{\ell_1^v}} .
\]
Recall that $ | {\mathbf a}-{\mathbf a(\mathbf p)} | _{\ell_1^v}=\sum
_{k=1}^p (1+k)^r |a_k-a_k(p)|+
\sum_{k=p+1}^{\infty}(1+k)^r |a_k|$\vspace*{2pt} in order to obtain from Lemma \ref
{lemma2} the desired assertion.
\end{pf*}
\begin{pf*}{Proof of Lemma \ref{lemma5}}
For simplicity, we write $p$ instead of $p(n)$.
From Lemma \ref{lemma3}, we have that the polynomial $A_p(z)$ has no
zeroes with magnitude less than or equal to
$1+1/p$. Since we easily get from assumption~(B) convergence of
$\widehat A_p(z)=1-\sum_{k=1}^p\widehat a_k(p)z^k$ to $A_p(z)$
uniformly\vspace*{1pt} on the closed disk with radius $1+1/p(n)$ the polynomial
$\widehat A_p(z)$ does not possess zeroes with magnitude less than or
equal to $1+1/p(n)$. Therefore, Cauchy's inequality for holomorphic
functions applies and yields
\begin{eqnarray*}
| \widehat\alpha_k(p)-\alpha_k(p) |
&\le& \frac{1}{(1+1/p)^k} \max_{|z|=1+1/p} | A_p(z)^{-1} -
\widehat A_p(z)^{-1} | \\
&= & \frac{1}{(1+1/p)^k} \max_{|z|=1+1/p} \frac{| A_p(z) -
\widehat A_p(z) | }{| A_p(z) \widehat A_p(z) | }
\\
&\le& \frac{1}{(1+1/p)^k} \max_{|z|=1+1/p} \frac{\sum
_{k=1}^p|\widehat a_k(p)-a_k(p)| (1+{1/p})^k }{| A_p(z)
\widehat A_p(z) | } \\
&=& \biggl( 1+\frac{1}{p}\biggr) ^{-k} \cdot\frac{1}{p^2}\cdot
{\mathcal O}_P(1) .
\end{eqnarray*}
\upqed\end{pf*}
\begin{pf*}{Proof of Theorem \ref{bootlinstat}}
A careful inspection of the proof of Theorem 3.3 in \citet
{Buhlmann97} [see also the corresponding
technical report \citet{Buhlmann95}] shows that only the following
properties of the underlying, the companion and the
fitted autoregressive process really are needed:
\begin{longlist}[(viii)]
\item[(i)] $\ve^{\ast}_t \xrightarrow{{\mathcal D^{\ast}}}
\widetilde\ve_t$ in probability,
\item[(ii)] $( X_{t_1}^{\ast},\ldots,X_{t_d}^{\ast})
\xrightarrow{{\mathcal D^{\ast}}} ( \widetilde X_{t_1},\ldots,
\widetilde X_{t_d}) $ in probability,\vspace*{2pt}
\item[(iii)] $\sum_{j=0}^{\infty} \vert\widehat\alpha_j(p(n))
- \alpha_j \vert\rightarrow_{n\to\infty} 0$ in probability,\vspace*{2pt}
\item[(iv)] $\sum_{j=0}^{\infty} j \vert\widehat\alpha
_j(p(n)) \vert$ is uniformly bounded in probability,\vspace*{2pt}
\item[(v)] $\sum_{j=0}^{\infty} j \vert\alpha_j \vert
< \infty$,
\item[(vi)] the empirical moments of $\widehat\ve_t(p(n))$ converge
for orders up to $2(h+2)$ to the moments of $\widetilde\ve_1$,
\item[(vii)] the autoregressive representation of infinite order of the process~$(\widetilde X_t)$ is invertible,
\item[(viii)] Yule--Walker parameter estimators are used for the
autoregressive fit of order $p(n)$ to the data $X_1,\ldots,X_n$.
\end{longlist}
Because of (A4) and (A5) and the easily obtained fact that
for the Mallows metric~$d_2$
%
%
\begin{equation}
d_2( \widehat F_n,F_n ) \rightarrow0 \qquad\mbox{in }
\mbox{ probability} ,
\end{equation}
where $\widehat F_n$ denotes the empirical distribution function of the
centered residuals $\widehat\ve_t(p(n))$,
$ t=p(n)+1,\ldots, n$ of the autoregressive fit and $F_n$ denotes the
empirical distribution function of fictitious observations $ \ve
_t(p(n)), t=p(n)+1,\ldots, n$, we obtain (i).

(ii) is obtained exactly along the lines as in Corollary 5.6 of
\citet{Buhlmann97}.

To see (iii), recall from Section \ref{generalAR} that we have
$\vert\widehat\alpha_j(p(n)) - \alpha_j(p(n))\vert=
(1+1/p(n))^{-j} \cdot1/p^2 \mathcal{O}_P(1) $ as well\vspace*{1pt} as
$\sum_{j=1}^{\infty} j \vert\alpha_j(p(n)) - \alpha_j\vert\le
C\cdot\sum_{j=p(n)}^{\infty} j \vert a _j \vert$ which converges
to $0$
as $n$ goes to infinity. These two assertions ensure (iii).

(iv) is
obtained exactly along these lines by using the fact that
$\sum_{j=1}^{\infty} j \vert\alpha_j\vert\le\infty$, cf.
(2.29), which also is (v).

Furthermore, it is easy to see that the difference of the empirical
moments (up to the necessary order) of
$\widehat\ve_t(p(n))$ and of $\ve_t$ converge to zero due to the
bounds for $\widehat\alpha_j(p)-\alpha_j(p)$,
$\alpha_j(p)-\alpha_j$ and $\alpha_j$, cf. (\ref{alphapdach}), (\ref
{alphap}) and (\ref{alpha}).
Together with (A5), we obtain (vi).

Finally, we use for the autoregressive fit Yule--Walker parameter
estimators and the autoregressive representation of
$(\widetilde X_t)$ is invertible (cf. Section~\ref{generalAR}). This
concludes the proof of Theorem \ref{bootlinstat}.
\end{pf*}
\begin{pf*}{Proof of Theorem \ref{intperiodoCLT}}
Also due to the results of \citet{Dahlhaus85}, we obtain that the
distribution of
$\sqrt{n}(M(\widetilde I_n,\varphi)-M(f_X,\varphi))$, where
$\widetilde I_n$ denotes the periodogram of $n$ observations of the
autoregressive companion process $(\widetilde X_t)$, cf. (\ref{ARcompanion}),
asymptotically is normal with mean zero and variance
%
%
\begin{equation}
\label{varianceintperiodotilde}
\bigl(E\ve_t^4/(E(\ve_t)^2)^2 -3\bigr) \biggl( \int_0^{\pi} \varphi(\lambda
)f_X(\lambda) \,d\lambda\biggr) ^2 + 2\pi\int_0^{\pi} \varphi
^2(\lambda)
f_X^2(\lambda) \,d\lambda.\hspace*{-28pt}
\end{equation}
Thus, it suffices to show that the distribution of the bootstrap
approximation of
$\sqrt{n}( M(\widetilde I_n^{\ast},\varphi)-M(\widehat f_{\mathrm
{AR}},\varphi)) $
shares the same asymptotic distribution.
Exactly along the lines of proof of Theorem 4.1 in \citet
{KrePap03} (without the additional nonparametric correction
considered therein) which makes use of Proposition 6.3.9 of \citet
{BrockwellDavis91}, we obtain the the desired result.
\end{pf*}
\begin{pf*}{Proof of Theorem \ref{intperiodoCLTspec}}
Let $ \widetilde{f}_n(\lambda)=\int_{-\pi}^{\pi} K_h(\omega-\lambda
)\widetilde{I}_n(\lambda) \,d\lambda$ and consider $ \widetilde
{Y}_n=\sqrt{nh}(\widetilde{f}_n(\lambda)-f_X(\lambda)) $. Since
$\widetilde{Y}_n$
converges to a Gaussian distribution
with mean and variance as in (\ref{biassd}) and (\ref{varsd}),
respectively, it suffices to show that $ \sqrt{nh}(f^\ast_n(\lambda
)-f_{\mathrm{AR}}(\lambda) )$ shares exactly the same asymptotic behavior
as $Y_n$. This however, follows exactly along the same lines as in the
proof of Theorem 5.1 in \citet{KrePap03}, again without the
additional nonparametric correction
considered therein.
\end{pf*}
\end{appendix}
%

\section*{Acknowledgments}
The authors would like to express their gratitude to two anonymous
referees. Their careful reading of a previous version of the paper and
their thorough reports led to a considerable improvement of the present
paper.


%

%
\printaddresses

\end{document}